\documentclass[12pt,psamsfonts]{amsart}
\usepackage{amsmath,amsthm,amsfonts,amssymb,verbatim}
\usepackage[numeric]{amsrefs}
\usepackage{eucal}
\usepackage{graphicx}
\usepackage[all,knot]{xy}
\xyoption{arc}

\newcommand{\CC}{\mathcal{C}}

\newcommand{\om}{\omega}

\newcommand{\Rep}{{\rm Rep}}

\newcommand{\N}{\mathbb{N}}
\newcommand{\Q}{\mathbb Q}
\newcommand{\R}{\mathbb{R}}

\DeclareMathOperator{\FPdim}{FPdim}

\DeclareMathOperator{\Hom}{Hom}

\DeclareMathOperator{\FSexp}{FSexp}
\newcommand{\one}{\mathbf{1}}

\newcommand{\Z}{\mathbb Z}

\newcommand{\ot}{\otimes}

\makeatletter

\newcommand{\Rmnum}[1]{\expandafter\@slowromancap\romannumeral #1@}
\makeatother

\numberwithin{equation}{section}

\newtheorem*{thm}{Theorem}

\newtheorem{theorem}{Theorem}[section]

\newtheorem{lemma}[theorem]{Lemma}

\newtheorem{question}[theorem]{Question}

\newtheorem{prop}[theorem]{Proposition}
\theoremstyle{definition}
\newtheorem{alg}[theorem]{Algorithm}

\newtheorem{remark}[theorem]{Remark}

\begin{document}
\title[]
{Modular categories, integrality and Egyptian fractions}
\author{Paul Bruillard}
\email{paul.bruillard@math.tamu.edu}
\address{Department of Mathematics\\
    Texas A\&M University \\
    College Station, TX 77843\\
    U.S.A.}
\author{Eric C. Rowell}
\email{rowell@math.tamu.edu}
\address{Department of Mathematics\\
    Texas A\&M University \\
    College Station, TX 77843\\
    U.S.A.}

\subjclass[2000]{Primary 18D10; Secondary 16T05, 11Y50}
\begin{abstract} It is a well-known result of Etingof, Nikshych and Ostrik that there are finitely many inequivalent integral modular categories of any fixed rank $n$.  This follows from a double-exponential bound on the maximal denominator in an Egyptian fraction representation of $1$.  A na\"ive computer search approach to the classification of rank $n$ integral modular categories using this bound quickly overwhelms the computer's memory (for $n\geq 7$).  We use a modified strategy: find general conditions on modular categories that imply integrality and study the classification problem in these limited settings. The first such condition is that the order of the twist matrix is $2,3,4$ or $6$ and we obtain a fairly complete description of these classes of modular categories.  The second condition is that the unit object is the only simple non-self-dual object, which is equivalent to odd-dimensionality.  In this case we obtain a (linear) improvement on the bounds and employ number-theoretic techniques to obtain a classification for rank at most $11$ for odd-dimensional modular categories.
\end{abstract}
\maketitle

\section{Introduction}
The problem of classifying low-rank modular categories has its roots in the classification problem for rational conformal field theories going back to the 1980s (predating the definition of modular category \cite{Tur92}).  Currently, the most complete results are in \cite{RSW} where unitary modular categories of rank at most $4$ are classified.  More generally, Ostrik classified ribbon fusion categories of rank at most $3$ \cite{Ostrik2} and fusion categories of rank at most $2$ \cite{Ostrik1}.  

By a generalized form of Ocneanu rigidity \cite[Prop. 2.31]{ENO} one may classify modular categories of a given rank up to finite ambiguity by classifying their Grothedieck semirings.  In \cite{HR} such an approach yielded a classification (up to Gr.-semirings) of modular categories of rank at most $5$ with the property that some  object is not isomorphic to its dual object.

Wang has conjectured that there are only finitely many inequivalent modular categories of each rank (see \cite[Conjecture 6.1]{RSW}), and the results mentioned above bear this out for rank at most $4$ and for rank at most $5$ in case some object is non-self-dual.  The most general class of modular categories for which Wang's conjecture has been verified is  for \emph{weakly integral} categories, that is, categories $\CC$ with $\FPdim(\CC)\in\N$ \cite[Prop. 8.38]{ENO}.  The proof relies upon a classical result of Landau \cite{Landau} that the diophantine equation:
\begin{equation}\label{sumeq}
1=\sum_{i=1}^{n}\frac{1}{x_{i}}
\end{equation}
has finitely many solutions with $x_i\in\N$.   Such solutions are \emph{Egyptian fraction} representations of $1$ which are of independent interest in combinatorial number theory (see \cite[Section D11]{Guy} and \cite[Seq. A002966]{OEIS}).  Moreover the number of solutions to (\ref{sumeq}) is at least exponential in $n$ since $\frac{1}{x}=\frac{1}{x+1}+\frac{1}{x^2+x}$ so that each (non-constant) $n$-term Egyptian fraction representation of $1$ leads to at least $2$ $(n+1)$-term representation.  In \cite{HR} a computational approach to classifying \emph{integral} modular categories (that is, with $\FPdim(X)\in\N$ for all objects $X$) of  rank $n$ is suggested using two facts:
\begin{enumerate}
 \item \cite[Lemma 1.2]{EG}: $\sqrt{\frac{\FPdim(\CC)}{\FPdim(X_i)}}\in\N$ for all simple $X_i$
\item \cite{Curtiss}: $\dim(\CC)\leq u_n$ where $u_n$ is inductively defined by $u_1:=1$ and $u_k:=u_{k-1}(u_{k-1}+1)$.
\end{enumerate}
As $u_n$ is double exponential in $n$, a direct search for solutions by computer quickly becomes infeasible.  

We attempt to circumvent this computational obstacle by finding interesting conditions on modular categories that imply integrality, and then use these conditions to simplify the classification problem for these classes of categories.  Specifically, we have:
\begin{thm}
Suppose $\CC$ is a modular category such that either:
\begin{enumerate}
 \item[(a)] the twist matrix $T$ satisfies $T^N=I$ for $N\in\{2,3,4,6\}$ or
\item[(b)] the only simple object in $\CC$ satisfying $X\cong X^*$ is the unit object $\one$
\end{enumerate}
then $\CC$ is integral.
\end{thm}
Statement (a) is proved below in Theorem \ref{twistint} and was inspired by Davydov who posed the question for $N=2$ to the second author.  Statement (b) is \cite[Theorem 2.2]{HR} and follows from a Galois theory argument.

We are interested in classifying categories with one of these two properties (a) or (b).  We obtain a fairly explicit description of modular categories with property (a) in Theorem \ref{N2346}.  Modular categories with property (b) of rank at most $11$ are shown to be \emph{pointed} in Theorem \ref{mnsdthm}, that is, $\FPdim(X_i)=1$ for each simple object $X_i$.  While property (b) may seem to be a rare condition at first glance one can show (see Prop. \ref{equivconds}) that it is equivalent to the condition that $\FPdim(\CC)$ is odd.

The remainder of the paper is structured as follows: in Section \ref{prelim} we collect together some notation and useful facts about modular categories.  We address the classification problem for modular categories having property (a) or (b) in Sections \ref{lowordertwist} and \ref{MNSD} respectively, and give some perspectives and futher directions in Section \ref{conclusions}.

\subsection*{Acknowlegements}
We benefitted from discussions with K. Rusek, M. Papanikolas, Y. Sommerh\"auser and D. Naidu.  The second author was partially supported by NSA grant H98230-10-1-0215.

\section{Notation}\label{prelim}
A modular category $\CC$ is a non-degenerate braided, balanced fusion category (see \cite{BK} or \cite{Rsurvey} for the complete axiomatic definition).  In this section we establish notation and describe some of the algebraic data and relations coming from the axioms of modular categories.  

We shall typically adopt the notation and normalizations of \cite{NS} for the data of a modular category.  The fusion coefficients are $N_{i,j}^k:=\dim\Hom(X_i\ot X_j,X_k)$ where $\one=X_0,\ldots,X_{n-1}$ are the (isomorphism classes of) simple objects the number of which ($n$) is called the \emph{rank} of $\CC$.  The diagonal twist matrix $T_{ij}:=\delta_{ij}\theta_i$ has finite order and the $S$-matrix is normalized so that $S_{00}=1$.  We will denote by $d_i$ the dimension of the simple object $X_i$, i.e. $d_i=S_{i0}=\dim(X_i)$.  Defining the \textit{Gauss sums} by $p_+=\sum_i \theta_i (d_i)^2$ and $p_-=\overline{p_+}$ we have $(ST)^3=p_+S^2$, $S^2=p_+p_-C$ where $C_{ij}=\delta_{ij^*}$ is the (involutive) charge conjugation matrix, which commutes with $T$.  In particular $(S,T)$ give rise to a (projective) representation of the modular group $SL(2,\Z)$.  We define the fusion matrices $(N_i)_{jk}:=N_{ik}^j$ and denote by $\FPdim(X_i)$ the largest eigenvalue of $N_i$, i.e. the \emph{FP-dimension} of $X_i$.  The global dimension of $\CC$ is $\dim(\CC)=\sum_i (d_i)^2$ and the global $FP$-dimension is $\FPdim(\CC):=\sum_i \FPdim(X_i)^2$.  If $\FPdim(\CC)=\dim(\CC)$ then $\CC$ is said to be \emph{pseudo-unitary}.  The \textit{Verlinde formula} relates the fusion coefficients to the $S$-matrix entries:
\begin{equation}\label{verlinde}
 N_{ij}^k=\sum_r\frac{S_{ir}S_{jr}\overline{S_{kr}}}{\dim(\CC)S_{0r}}.
\end{equation}
A pair of matrices $(S,T)$ satisfying the above relations such that the right-hand-side of (\ref{verlinde}) is a non-negative integer for all $i,j,k$ is called a \emph{modular datum} (\cite{Gan05}).  A modular category $\CC$ with corresponding $S$ and $T$ matrices is called a \emph{categorification} of $(S,T)$.

A category $\CC$ is called \emph{integral} if $\FPdim(X_i)\in\N$ for all simple objects $X_i$.  All categories encountered in this work will be integral and hence pseudo-unitary so that we may assume that $\FPdim(X_i)=d_i$ and $\FPdim(\CC)=\dim(\CC)$ by \cite[Prop. 8.23, 8.24]{ENO}.

The entries of $N_i$, $S$ and $T$ satisfy further relations:
\begin{equation}\label{GS}
 p_+p_-=\dim(\CC)
\end{equation}
and
\begin{equation}\label{twisteq}
  \theta_i\theta_jS_{ij}=\sum_k N_{i^*j}^k\theta_k d_k.
\end{equation}
Since $\dim$ is a character of the Grothendieck semiring of $\CC$ we have:

\begin{equation}\label{chars}
\sum_k N_{ij}^k d_k=d_id_j.
\end{equation}

The (second) \textit{FS-indicator} is defined to be
 $$\nu_k=\frac{1}{\dim(\CC)}\sum_{i,j}N_{i,j}^kd_id_j\left(\frac{\theta_i}{\theta_j}\right)^2$$
 and satisfies:
\begin{equation}\label{fseq}
 \nu_k=\begin{cases} 0 & X_k\not\cong X_k^*\\ \pm 1 & X_k\cong X_k^*\end{cases}.
\end{equation}

More generally the $n$-th FS-indiactor $\nu_k^{(n)}$ may be defined for simple objects in modular categories in an analogous way (\cite[Theorem 7.5]{NS1}).  The quantity $\FSexp(\CC)$ is the smallest integer $m>0$ such that $\nu_k^{(m)}=d_k$ for all simple $k$.  For modular categories $\FSexp(\CC)$ coincides with the order of the $T$ matrix (\cite[Theorem 7.7]{NS1}).

\section{Low-order Twist Matrices}\label{lowordertwist}
In this section we study modular categories with twist matrix of order $2,3,4$ or $6$ of arbitrary rank.  Our first result is:
\begin{theorem}\label{twistint}
Suppose that $T$ is the twist matrix of a modular category $\CC$ such that $T^N=I$ for $N\in\{2,3,4,6\}$.  Then $\CC$ is integral.
\end{theorem}
\begin{proof} First we observe that by \cite[Prop. 5.7]{NS} the entries of the $S$-matrix for $\CC$ must lie in $\Q\left(\theta_1,\ldots,\theta_{n-1}\right)=\Q\left(\zeta_N\right)$ where $\zeta_N$ is a primitive $N$th root of unity.  Since $s_{ij}$ are algebraic integers and $\varphi(N)\leq 2$ (Euler's totient $\varphi$), if $s_{ij}\in\R$ then $s_{ij}\in\Z$.  In particular some column of the $S$-matrix must be an integer multiple of the vector of $FP$-dimensions (since $\CC$ is modular) and so $\FPdim(X_i)$ are rational integers.
\end{proof}

We now characterize modular categories with $T^N=I$ for $N \in \{2,3,4,6\}$ in the following:
\begin{theorem}\label{N2346} Suppose $\CC$ is a modular category such that $T^N=I$.  Then: 
\begin{enumerate}
\item[(a)] If $N=2$ then $\CC$ is braided tensor equivalent to a subcategory of $\Rep(D^\om G)$ where $G$ is an abelian $2$-group of exponent $2$ and $\dim(\CC)=2^{2s}$ (in particular, $\CC$ is pointed),
\item[(b)] If $N=3$ then $\CC$ is braided tensor equivalent to a subcategory of $\Rep(D^\om G)$ where $G$ is a $3$-group of exponent $3$, and 
    \item[(c)] If $N=4$ then $\CC$ is braided tensor equivalent to a subcategory of $\Rep(D^\om G)$ where $G$ is a $2$-group of exponent $2$ or $4$.
        \item[(d)] If $N=6$ then $\CC$ is solvable and hence weakly integral.
\end{enumerate}
\end{theorem}
\begin{proof}

If  $T^2=I$ then each $X_k$ is self-dual since the entries of the $S$-matrix are real.   Computing the FS-indicator we have:
$$\nu_k=\frac{1}{\dim{\CC}}\sum_{i,j} N_{ij}^kd_id_j\left(\frac{\theta_i}{\theta_j}\right)^2=\pm 1$$
for all $k$.  But $\left(\frac{\theta_i}{\theta_j}\right)^2=1$ by assumption and $\sum_iN_{ij}^kd_i=d_kd_j$ (as $N_{ij}^k$ is totally symmetric in the self-dual case) so we may simplify:

$$\nu_k=\frac{1}{\dim{\CC}}\sum_{j} d_kd_j^2=d_k.$$
But $\nu_k=\pm 1$ and $d_k>0$ so this implies $d_k=1$ for all simple objects $X_k$.  We must then have $X_k^{\ot 2}\cong \one$ as well.  Thus $\CC$ is pointed and has the same fusion rules as the group $\Z_2^m$ where $2^m$ is both the rank and global dimension of $\CC$.  Now since $p_+\in \R$ we have $p_+=p_-\in\Z$ so that by (\ref{GS}) we must have $\dim(\CC)=2^m=(p_+)^2$ and hence $m=2s$.  Since $\CC$ is braided tensor equivalent to a subcategory of $Z(\CC)\cong\Rep(D^\om G)$ we must have $G$ an elementary abelian $2$-group.

For (b) and (c) with $N=3$ (resp. $4$) we use \cite[Theorem 8.4]{NS1} to conclude that $\dim(\CC)=p^t$ where $p=3$ (resp. $2$).  In particular $\CC$ is braided and nilpotent and hence group theoretical by \cite[Theorem 6.10]{DGNO1}. Thus we conclude that $\CC$ is braided tensor equivalent to a subcategory of $Z(\CC)\cong\Rep(D^\om G)$ where $G$ is a $p$-group and since $\FSexp(\CC)=\FSexp(Z(\CC))=3$ (resp. $4$) by \cite[Corollary 7.8]{NS1} $G$ must have exponent $3$ (resp. $2$ or $4$ as the exponent of $G$ must divide that of $\Rep(D^\om G)$).

For (d) with $N=6$ \cite[Theorem 8.4]{NS1} implies that $\dim(\CC)=2^a3^b$ with $a,b>0$ and so by \cite[Theorem 1.6]{ENO2} $\CC$ is solvable, and in particular weakly group-theoretical.
\end{proof}

\begin{remark}\ \\
\begin{itemize}
 \item In principle, Theorem \ref{N2346}(a),(b),(c) can be used to completely classify modular categories with twist matrix of orders $2,3$ and $4$.
\item Theorem \ref{N2346}(d) is sharp in the sense that we cannot conclude $\CC$ is group theoretical as \cite[Example 5.3(a)]{GNN} shows that non-group-theoretical modular categories with $T^6=I$ exist.
\item In \cite[Prop. 5.1]{Cuntz} an infinite family of modular data $(S,T)$ depending on an integer $k\geq 0$ is described with $S_{ij}\in\Z$ for all $i,j$ and $T^2=I$.  Theorem \ref{N2346}(a) shows that this family is not categorifiable for $k\geq 1$.
\end{itemize}
\end{remark}

\section{Integrality and Egyptian Fractions}\label{MNSD}
In this section we classify integral modular categories of rank at most $6$ and maximally non-self-dual modular categories of rank at most $11$.

If $\CC$ is a (pseudo-unitary) integral modular category with dimensions of simple objects $d_1,d_2,\ldots,d_{n}=\dim(\one)=1$ then by \cite[Lemma 1.2]{EG} the numbers $x_i=\dim(\CC)/(d_i)^2$ are \emph{integers} and so satisfy eqn. (\ref{sumeq}).  In \cite[Lemma 2.3]{HR} it is shown that after relabelling if necessary so that $x_1\leq x_2\leq\cdots\leq x_n=\dim(\CC)$ then:
$$k\leq x_k\leq u_k(n-i+1)\quad \mathrm{for all}\quad 1\leq k\leq n,$$
where $u_1=1$, $u_k=u_{k-1}(u_{k-1}+1)$.
These bounds rely upon two classical results in \cite{Landau} and \cite{Curtiss}.  Notice also that $x_{n}/x_i=(d_i)^2$ is a perfect square for each $i$.  The \emph{trivial solution} $x_1=\cdots=x_n=n$ corresponds to a \emph{pointed} modular category of rank $n$ as then $d_i=1$ for all $i$.  Our basic algorithm for finding $x_1\leq\cdots\leq x_n$ satisfying these conditions is as follows:
\begin{alg}\label{algorithm1} \ \\
\begin{enumerate}
\item form the set of pairs $S_2:=\{(j,j\cdot i^2):2\leq j\leq n, 1\leq i\leq \sqrt{u_n/j}\}$.  These are the possible pairs $(x_1,x_n)$.
\item having determined $S_k$ for some $2\leq k\leq n-2$ let $S_{k+1}$ be the set of $(k+1)$-tuples $(j_1,\ldots,j_{k-1},J,j_{k})$ such that 
\begin{enumerate}
\item $(j_1,\ldots,j_k)\in S_k$,
\item $\max(j_{k-1},k)\leq J\leq \min(j_k,u_k(n-k+1))$,
\item $\sqrt{j_k/J}\in\Z$ and
\item $\frac{1}{J}+\sum_{i=1}^k \frac{1}{x_i}\leq 1$.
\end{enumerate}
\item the solution set $S_n$ is the set of $n$-tuples $(j_1,\ldots,j_{n-2},J,j_{n-1})$ where 
\begin{enumerate}
\item $(j_1,\ldots,j_{n-1})\in S_{n-1}$,
\item $\max(j_{n-2},n-1)\leq J\leq \min(j_{n-1},2u_{n-1})$,
\item $\sqrt{j_{n-1}/J}\in\Z$ and 
\item $\frac{1}{J}+\sum_{i=1}^{n-1} \frac{1}{x_i}=1$.
\end{enumerate}
\end{enumerate}
\end{alg}

For $n=5$ and $n=6$ the solutions are: $\{(5,5,5,5,5),(2,8,8,8,8)\}$ and $$\{(6,6,6,6,6,6),(3,3,12,12,12,12),(4,4,4,9,9,36)\}$$ respectively.  These computations were carried out in a few hours on a desktop computer using Maple, and we have:
\begin{theorem}
Any integral modular category of rank at most $6$ is pointed.
\end{theorem}
\begin{proof}
For rank at most $5$ this result is known, see: \cite{RSW} and \cite{HR}. For rank $6$ we need only show that no modular category may have simple dimensions $(1,1,1,1,2,2)$ or $(1,2,2,3,3,3)$ corresponding to the non-pointed solutions above.  The first case is eliminated by \cite[Prop. 4.11]{propF}: integral modular categories $\CC$ with $\FPdim(\CC)=2^2 3$ must be pointed.  On the other hand, a fusion category of dimension $36$ must be solvable and contain a non-trivial invertible object by \cite[Prop. 8.2]{ENO2} which eliminates the second solution.
\end{proof}
We remark that there is a non-pointed integral modular category of rank $8$, namely $\Rep(D(S_3))$, the representation category of the double of the symmetric group $S_3$.  We expect that this is the smallest possible rank for non-pointed integral modular categories.
Unfortunately for $n=7$ one has $u_7\approx 10^{13}$ and attempts to implement Algorithm \ref{algorithm1} quickly overwhelms a computer's memory. 
This motivates passing to a restricted class of modular categories that enjoy integrality.

With a view towards improving the upper bound $u_n$ above, suppose that $\CC$ is a maximally non-self dual (MNSD) modular category, i.e. each non-trivial simple object $X_i$ satisfies $X_i\not\cong X_i^*$.  Then $\CC$ is integral with rank $n=2k+1$ by \cite[Theorem 2.2]{HR}.  Since $d_i=\dim(X_i)=\dim(X_i^*)$ we may label the dimensions of simple objects: $d_1,d_1\ldots,d_k,d_k,d_{k+1}=\dim(\one)=1$ so that $\FPdim(\CC)=1+2\sum_{i=1}^k d_i^2$ is odd.  After reordering the $k+1$ integers $x_i:=\dim(\CC)/d_i^2$ so that $x_1\leq x_2\leq\cdots\leq x_k\leq x_{k+1}=\dim(\CC)$ they satisfy the special case of eqn. (\ref{sumeq}):
\begin{equation}\label{newsumeq}
\frac{1}{x_{k+1}}+\sum_{i=1}^{k}\frac{2}{x_k}=1
\end{equation}
The maximal non-self-dual condition may seem somewhat exceptional, but it follows from  \cite[Theorem 2.2]{HR}, \cite[Corollary 8.2(ii)]{NS1} and \cite[Corollary 2.33]{DGNO2} that:
\begin{prop}\label{equivconds} Let $\CC$ be a modular category of rank $n$.  The the following are equivalent:
\begin{enumerate}
 \item $\CC$ is maximally non-self-dual
\item  $\FPdim(\CC)$ is odd
\item $\CC$ is (monoidally) equivalent to $\Rep(A)$ where $A$ is an odd-dimensional, semisimple quasi-Hopf algebra.
\end{enumerate}
\end{prop}

The following generalization of \cite[Lemma 2.3]{HR} gives a linear improvement of the upper bounds in Algorithm \ref{algorithm1} for a restricted class of integral modular categories and includes the MNSD setting as the special case $\ell=2$:

\begin{lemma}
\label{new bounds 1}
Suppose that $\mathcal{C}$ is an integral modular category of rank $n$ and there exist $k$ (weakly decreasing) integers
$p_{1}\geq p_{2}\geq \cdots\geq p_{k}$ such that the non-trivial (isomorphisms classes of) simple objects can be partitioned into $k$ sets $P_1,\ldots,P_k$ such that $X\in P_i$ has $\dim(X)=p_i$ and $|P_i|=\ell$ for all $1\leq i\leq k$.
 Then:
\begin{itemize}
\item[(i)] The numbers $x_{i}:=\frac{\dim{\mathcal{C}}}{p_{i}^{2}}$ form a weakly increasing sequence of integers such that $\sum_{i=1}^{k}\frac{\ell}{x_{i}}+\frac{1}{x_{k+1}}=1$, and
\item[(ii)] $\ell i\leq x_{i}\leq\left(n+\ell-i\ell\right)\frac{A_{i}}{\ell}$ for $i\leq k$ and $(k+1)\leq x_{k+1}\leq\frac{A_{k+1}}{\ell}$ where $A_{1}=\ell$ and $A_{i}=A_{i-1}\left(A_{i-1}+1\right)$.
\end{itemize}
\end{lemma}
\begin{proof}
The proof of (i) proceeds exactly as in \cite[Lemma 2.3]{HR} so we focus on (ii).  First define $x_i^\prime=x_i$ for $1\leq i\leq k$ and $x_{k+1}^\prime=x_{k+1}/\ell$ so that:
\begin{equation}\label{sumwithell}
\sum_{i=1}^{k+1}\frac{1}{x_i^\prime}=\frac{1}{\ell}
\end{equation}
and $x_1^\prime\leq\cdots\leq x_{k+1}^\prime=\ell x_{k+1}$.  Observing that $n=k\ell+1$ we finds that (ii) becomes $\ell i\leq x_{i}^\prime\leq\left(n+\ell-i\ell\right)\frac{A_{i}}{\ell}$.
Since the $x_i^\prime$ are weakly increasing the lower bounds are clear: for if $\ell i>x_i^\prime$ for some $i$ then $\ell i> x_i^\prime\geq x_{i-1}^\prime\geq\cdots\geq x_1^\prime$ so $1/\ell=i\frac{1}{\ell i}>\sum_{j=1}^{i}\frac{1}{x_j^\prime}$, a contradiction.  For the upper bounds follow the same strategy as in \cite{HR}: we use Takenouchi's \cite{Ta} bound on the largest denominator $x_{k+1}^\prime$ in eqn. (\ref{sumwithell}) and Landau's \cite{Landau} estimate.  In particular, if
$$\sum_{i=1}^k\frac{1}{y_i}=r$$ with $y_1\leq y_2\leq\cdots\leq y_k$ then Landau's result says $y_i\leq (k-i+1)/r_{i-1}$ where $r_0=r$ and $r_i=r_{i-1}-1/y_i$.  Takenouchi's result says that the maximum denominator $y_i$
of a solution to $\frac{1}{\ell}=\sum_{i=1}^k\frac{1}{y_i}$ is $A_k$ where $A_1=\ell$ and $A_i:=A_{i-1}(A_{i-1}+1)$.  So in particular $1/r_{i-1}\leq A_{i}$ for all $i$ and the result follows.

\end{proof}

Lemma \ref{new bounds 1} provides a linear improvement on the bounds in \cite[Lemma 2.3]{HR}: one has roughly $\mathrm{rank}(\CC)/\ell$ variables as opposed to $\mathrm{rank}(\CC)$.
Taking the case $\ell=2$ one can modify Algorithm \ref{algorithm1} to search for possible MNSD modular categories in the obvious way.  As the bounds are still double exponential the best result we can achieve with our (modified) algorithm is:     
\begin{theorem}\label{mnsdthm}
All maximally non-self dual modular categories of rank at most $11$ are pointed.
\end{theorem}
\begin{proof}  The only non-trivial solution the modified algorithm produces is in rank $11$, consisting of $9$ simple objects of dimension $1$ and two simple objects of dimension $3$.  To eliminate this possibility we use a result in \cite{GN}: any integral modular category $\CC$ is faithfully graded by its universal grading group $U(\CC)$, which is isomorphic to the group of invertible objects in $\CC$.  Since each component of a faithful grading must have the same $FP$-dimension we see that this is impossible as $|U(\CC)|=9$ and any component $\CC_g$ of the grading containing a simple object $X$ with $\FPdim(X)=3$ must have $\FPdim(\CC_g)\geq 9$.
\end{proof}

\begin{remark}
 An obvious source of MNSD modular categories are those of the form $\Rep(D^\om G)$ where $|G|$ is odd (that is, the representation category of the twisted double of a finite group of odd order).  
The smallest rank example of a non-pointed MNSD modular category we are aware of is the rank $25$ category $\Rep(D(\Z_3\ltimes\Z_7))$ of dimension $441$ where $\Z_3\ltimes\Z_7$ is a non-abelian semidirect product.  One might speculate that any MNSD modular category of rank at most $23$ is pointed. 
\end{remark}

We close with the following related question that we find interesting:
\begin{question}
 Is there an odd dimensional modular category that is non-group-theoretical?
\end{question}

\section{Conclusions and Future Directions}\label{conclusions}
We have made significant progress towards classifiying modular categories with twist matrices of order $2,3,4$ or $6$.  It is not clear how to extend these techniques to other orders.  One might hope that one could generalize \cite[Theorem 8.4]{NS1} to non-integral categories to get some statement about the possible prime divisors of $\FPdim(\CC)$ in terms of those of the order of the twist matrix.  However, the relevant number field would no longer be $\Z$ so that other complications would arise particularly concerning units.

For rank greater than $11$ even the improved bounds for MNSD modular categories become too large for our computational techniques.  However, if one considers MNSD modular categories with:
\begin{enumerate}
 \item rank $r$
\item  $\FPdim(X_i)\leq M$, for each simple $X_i$ and
\item $|\{\FPdim(X_i):0\leq i\leq r-1\}|=t$  ($X_i$ simple).
\end{enumerate}
further headway can be made.  For example, taking $r\leq 23$, $M=45$ and $t\leq 6$ we set $k_i=|\{j:\FPdim(X_j)=d_i\}|$ with $d_0=1$ we obtain (using \cite[Lemma 1.2]{EG} and \cite[Prop. 8.15]{ENO}):
\begin{eqnarray*}
 &&k_0\equiv 1\mod{2}\\
&&k_i\equiv 0\mod{2}, (i\neq 0)\\
&&\sum_{i=0}^{t-1} k_i=r\\
&&\sum_{j\neq i} k_j(d_j)^2\equiv0\mod{(d_i)^2}\\
&&\sum_{j\neq 0} k_j(d_j)^2\equiv0 \mod{k_0}
\end{eqnarray*}
For any fixed $r\leq 23$ and subsets of $\{d: 3\leq d\leq 45,\/ d\in 2\Z+1\}$  of size at most $5$
we convert these congruences to a system of linear diophantine equations and then use the Smith normal form of the corresponding matrix to solve for the $k_i$.  Applying 
 further classification theorems and \emph{ad hoc} techniques produces the following partial result:
\begin{theorem}
If $\mathcal{C}$ is a non-pointed maximally non-self dual modular category with $13\leq\text{rank}\left(\mathcal{C}\right)\leq23$ then either $\CC$ has:
\begin{enumerate}
\item[(a)] $\FPdim(X_i)\geq 47$ for some simple $X_i$ or
\item[(b)]  $|\{\FPdim(X_i):0\leq i\leq m-1\}|\geq 7$
\end{enumerate}
\end{theorem}
This is further evidence that perhaps a non-pointed MNSD modular category of smallest rank is indeed the category of rank $25$ and dimension $441$ described above.

More generally, one can also obtain partial classification results for integral modular categories of small rank by bounding the $FP$-dimension further (below the double-exponential bound of $u_n$ in Algorithm \ref{algorithm1}).  For example, using a streamlined version of Algorithm \ref{algorithm1} provided by K. Rusek and some classification theorems we can prove:
\begin{theorem}
 If a non-pointed modular category $\CC$ has rank $7$ then $\FPdim(\CC)>10^{5}$.
\end{theorem}

\bibliographystyle{amsxport}
\bibliography{Refs}

\end{document}